\documentclass[aos,MSNbibl,nameyear,dvips]{arximspdf}
\usepackage{graphicx}
\doi{10.1214/13-AOS1186} 
\volume{42}
\issue{1}
\pubyear{2014}
\firstpage{190}
\lastpage{210}

\makeatletter
\newcommand{\rrVert}{\Vert}
\newcommand{\rrvert}{\vert}
\newcommand{\llVert}{\Vert}
\newcommand{\llvert}{\vert}
\newtheorem{theorem}{Theorem}[section]
\newtheorem{lemma}{Lemma}[section]
\newproclaim{remark}{Remark}[section]
\makeatother

\begin{document}
\begin{frontmatter}

\title{Statistical inference based on robust low-rank data~matrix
approximation}
\runtitle{Robust approximations with low-rank data matrices}

\begin{aug}
\author[A]{\fnms{Xingdong} \snm{Feng}\thanksref{m1}\ead[label=e1]{feng.xingdong@mail.shufe.edu.cn}}
\and
\author[B]{\fnms{Xuming} \snm{He}\corref{}\thanksref{m2}\ead[label=e2]{xmhe@umich.edu}}
\runauthor{X. Feng and X. He}
\affiliation{Shanghai University of Finance and Economics and
University of Michigan}
\address[A]{School of Statistics and Management\\
Shanghai University of Finance and Economics\\
and Key Laboratory of Mathematical Economics (SUFE)\\
Ministry of Education\\
777 Guoding Road\\
Shanghai 200433\\
China\\
\printead{e1}} 
\address[B]{Department of Statistics\\
University of Michigan\\
439 West Hall\\
1085 South University Avenue\\
Ann Arbor, Michigan 48109\\
USA\\
\printead{e2}}
\end{aug}
\thankstext{m1}{Supported by NNSF of China Grant 11101254, Shanghai
Pujiang Program, Program for Changjiang Scholars and Innovative Research
Team in Universities, and Key Laboratory of Mathematical Economics (SUFE),
Ministry of Education, China.}
\thankstext{m2}{Supported by NSF Grants DMS-13-07566, DMS-12-37234, NIH
Grant R01GM080503 and NNSF of China Grant 11129101.}

\received{\smonth{8} \syear{2013}}

%
\begin{abstract}
The singular value decomposition is widely used to approximate data
matrices with lower rank matrices.
Feng and He [\textit{Ann. Appl. Stat.} \textbf{3} (2009) 1634--1654]
developed tests on dimensionality of the mean structure of
a data matrix based on the singular value decomposition.
However, the first singular values and vectors can be driven by a small
number of outlying measurements.
In this paper, we consider a robust alternative that moderates the
effect of outliers in low-rank approximations. Under
the assumption of random row effects, we provide the asymptotic
representations of the robust low-rank approximation.
These representations may be used in testing the adequacy of a low-rank
approximation. We use oligonucleotide gene microarray
data to demonstrate how robust singular value decomposition compares
with the its traditional counterparts. Examples show that the robust
methods often lead to a more meaningful assessment of the
dimensionality of gene intensity data matrices.
\end{abstract}

%

\end{frontmatter}

\section{Introduction}\label{intro}
Research on robustness dates back to the prehistory of
statistics. However, the concepts and theories of robust statistics
have not
been formally and systematically established until recent decades
[\citet{Huber}, \citet{Hampel}]. Much work on robust
statistics has focused on
linear regression and multivariate location-scatter models. It has been
well recognized that the least squares method under those models is
sensitive to a small number of outliers. Robust methods are generally
developed to down-weight outliers.

The singular value decomposition (SVD) of a~data matrix is often used
as a~data reduction tool.
In fact, the SVD can be viewed as a basic tool in dimension reduction.
Consider a data matrix
\[
\mathbf{Y}=\pmatrix{ y_{11}&\cdots&y_{1m}
\cr
\vdots&\vdots&
\vdots
\cr
y_{n1}&\cdots&y_{nm}},
\]
of $n$ rows and $m$ columns, where $m$ is fixed. An approximation of
rank $r$ to the matrix can be found by
%
\begin{equation}
\label{svd1} \sum_{i=1}^n\sum_{j=1}^m
(y_{ij}- z_{ij})^2,
\end{equation}
where $z_{ij}$ are the elements of $\mathbf{Z} = R C$ for an $n \times
r$ matrix $R$ and a $r \times m$ matrix~$C$. The matrices $R$ and $C$
are not identifiable in this formulation, so additional constraints may
be imposed to ensure identifiability.
As pointed out in \citet{ammann} and \citet{ChenAndHe},
the SVD is equivalent to the least squares approach to a bilinear
regression model, so it suffers from the usual
lack of robustness against outliers. 

\citet{ruppert} have used the trimmed least squares estimation in the
linear model by using weights obtained from some initial consistent
estimates, and \citet{gervini} have considered a variant of the trimmed
method leading to the maximum breakdown point and full asymptotic
efficiency under normal errors. In this paper, we adopt the idea of
using trimmed least squares estimation, where the scheme of choosing
weights is explained in Section~\ref{init}. The low-rank approximation
of matrices by weighted least squares has been considered by \citet{gabriel}, but their weights are fixed, while the weights of the
proposed method in this paper are obtained from an initial robust estimate.


We will consider a two-step approximation method in this paper. More
specifically, we consider the first approximation by minimizing
%
\begin{equation}
\label{robobj1} \sum_{i=1}^n\hat{w}_i
\sum_{j=1}^m \Biggl(y_{ij}-\sum
_{k=1}^r\theta_{ki}
\phi_{kj} \Biggr)^2,
\end{equation}
where $\hat{w}_i$ are the weights based on an initial estimate (to be
described later), $\theta_{ki}$~are the elements of $R$,
and $\phi_{kj}$ are the elements of $C$.
However, it is clear that the estimates of $\theta$'s are the linear
combination of vectors $\underline{y}$'s given $\phi$'s, so it implies
that this lower-rank approximation is not robust against outliers.
Then we consider the second approximation by using the estimated $\phi
$'s from the first step, denoted collectively as $\underline{\tilde
{\phi}}_k$ ($k=1, \ldots, r$), and then minimizing
%
\begin{equation}
\label{robobj2} \sum_{i=1}^n\sum_{j=1}^mL
\Biggl(y_{ij}-\sum_{k=1}^r\theta
_{ki}\tilde\phi_{kj} \Biggr),
\end{equation}
over the $\theta$'s for some robust loss function $L$, where $\tilde
\phi
_{kj}$ is the $j$th component of~$\underline{\tilde{\phi}}_k$.
Our statistical analysis will be performed under
the following model:
%
\begin{equation}
\underline{y}_i=\sum_{k=1}^r
\theta_{ki}^{(0)}\underline{\phi}_k^{(0)}
+ \underline{\varepsilon}_i,\qquad i=1,\ldots,n, \label{equ0}
\end{equation}
where $\underline{y}_i=(y_{i1},\ldots,y_{im})^T$ is the $i$th
observed vector,
$\underline{\theta}_k^{(0)}=(\theta_{k1}^{(0)},\ldots,\theta_{kn}^{(0)})^T$
is used to explain the row effects, and
$\underline{\phi}_k^{(0)}=(\phi_{k1}^{(0)},\ldots,\phi_{km}^{(0)})^T$
is used to explain the column effects in the data matrix. The row
effects $\theta^{(0)}$'s are assumed to be random, and the length of
observed vectors, $m$, is fixed. We are interested in the structure of
the mean matrix $E(\mathbf{Y})$, and the uniqueness of the low-rank
representation is implied by conditions (M1) and (M2) given in Appendix~\ref{AssMod}. In our work, we assume that each component of
$\underline
{\varepsilon}_i$ in the model is symmetrically distributed, but
outliers might be present in the data. The robust methods are meant to
be reliable against violations of the model assumptions.

Our model includes that of \citet{feng} as a special case where
$\underline{\varepsilon}_i$ is Gaussian. The main contribution of the
present paper is to develop a robust procedure that can accommodate
outlying measurements in the data matrix. To achieve this goal, we have
to utilize nonlinear operations in the estimation procedure, and
consequently, we need to analyze the statistical properties of the
robust procedures with a new set of techniques.

When the data matrix is the sum of low-rank and sparse matrices, the
theory of the exact recovery of both matrices has been established by
\citet{candes} and \citet{zhou}. \citet{agarwal} further consider a
broader class of models, where random errors are introduced, and the
penalized method is used for estimation.
These authors have provided deterministic error bounds for their
estimates while allowing the number of columns to grow with~$n$, but in
this paper we are interested in hypothesis testing based on the
asymptotic representation of the robust estimates with a fixed number
of columns.

For the estimates of $\underline{\theta}_k^{(0)} $ and $\underline
{\phi
}_k^{(0)}$ ($k=1,\ldots, r$) obtained from (\ref{robobj2}) and (\ref
{robobj1}), respectively,
we shall derive their asymptotic representations
in Section~\ref{estimates} as $n \to\infty$.
In Appendix~\ref{app}, we discuss some finite sample properties of the
estimators (\ref{robobj2}), which are critical for the theoretical
development in Section~\ref{appli}, where we robustify the tests of
unidimensionality for testing the adequacy of a unidimensional model
against the alternative rank-two mean structure for the data matrix. 
In Section~\ref{application}, we compare the results of testing
unidimensionality of matrices from \citet{feng}
with those from the robust alternative in microarray data analysis.
Technical assumptions of our model are given in Appendix~\ref{ass}, and
the proofs for the lemmas and the theorems given in the paper can be
found in Appendix~\ref{app} and in the supplementary material
[\citet{feng1}].

\section{Estimation procedure}\label{init}
In this paper, we propose the following procedure to estimate the row
and column parameters of model (\ref{equ0}).

\begin{longlist}[Step 0.]
\item[Step 0.] Construction of an initial robust estimate of column parameters:

\begin{enumerate}[(I4)]
\item[(I1)]prechoose a constant $\alpha^* $ (typically between 0.1
and 0.5);
\item[(I2)]select $ \lceil(1-\alpha^*)n \rceil$ rows randomly
from the data matrix, and denote this matrix as $\mathbf{Y}^*$, where
$ \lceil x  \rceil$ is the smallest integer greater than $x$;
\item[(I3)]carry out the regular SVD on the matrix $\mathbf{Y}^*$ and
obtain the first $r$ right singular vector as $\underline{\hat{\phi
}}_k$, $k=1,\ldots,r$;
\item[(I4)] estimate the row parameters $\underline{\theta}_k^{(0)} $
by minimizing the objective function (\ref{robobj2}), in which
$\underline{\tilde{\phi}}_k$ is\vspace*{1pt} replaced by $\underline{\hat{\phi}}_k$.
The resulting estimate is denoted as~$\underline{{\hat{\theta}}}_k$;
\item[(I5)] repeat (I2)--(I4) for a prespecified number of times (to be
discussed later), and find the subset of $ \lceil(1-\alpha
^*)n
\rceil$ rows that gives the minimum value of $\sum _{i=1}^n\sum_{j=1}^m L  (y_{ij}-\sum_{k=1}^r\hat\theta_{ki}\hat\phi
_{kj} )$.
\end{enumerate}
\item[Step 1.] Computation of the weighted least squares to improve
efficiency of the column parameters:

\begin{enumerate}[(1a)]
\item[(1a)] given the initial estimate of the column parameters, choose
a trimming proportion $\alpha\le\alpha^*$ and calculate the weights
%
\begin{equation}
\label{weight} \hat{w}_i=1 \bigl(\hat\xi_\alpha<\|
\underline{ \hat{e}}_i\|^2\leq\hat \xi _{1-\alpha}
\bigr),
\end{equation}
where
$\hat{\xi}_{\alpha}$ is the sample $\alpha$ quantile of
$\|\underline
{\hat{e}}_i\|^2$ and $\underline{\hat{e}}_i= (I_m-\sum_{k=1}^r\underline{\hat{\phi}}_k\underline{\hat{\phi}}{}_k^T
)\underline{y}_i$;
\item[(1b)] given the weights, obtain the estimate $\tilde{\underline
{\phi}}_k$ ($k=1,\ldots,r$) of the column parameters by minimizing
(\ref
{robobj1}) over the row and column parameters.
\end{enumerate}
\item[Step 2.] Updating\vspace*{-1pt} row effect estimates with robustness: given
$\tilde{\underline{\phi}}_k$ from step~1, obtain the estimate
$\tilde{
\underline{\theta}}_{k}$ of the row effects by minimizing (\ref
{robobj2}) over the row parameters.
\end{longlist}

In step~0, we obtain an initial root-$n$ robust estimate of the column
parameters~$\underline{\phi}_k^{(0)}$, denoted as ${\underline{\hat
\phi
}}_k$, $k=1,\ldots,r$.
The choice of $\alpha^*$ should reflect what percentage of outlying
rows we expect, and it is similar to the amount of trimming one chooses
to use in the trimmed mean. The number of subsets used in (I5) is fixed
and should be chosen to ensure that there is a high probability that
one of the subsets contains no outliers. For example, if we have 20
rows in the data matrix and expect 2 outlying rows, by choosing $\alpha
^* =0.3$ to use subsets of 14 rows, the probability that one random
subset is outlier-free is nearly 0.08. If we use 100 random subsets in
(I5), the probability of having at least one outlier-free subset is
greater than 0.999. Simple calculations like this show that we can
obtain a robust estimate through this procedure with high probability.

Because the estimate $\underline{\hat{\phi}}_k$ in (I3) is the least
squares estimate considered in \citet{feng}, and the size of the subset
is proportional to $n$, then the initial estimate of column vectors
here is root-$n$ consistent. Given the initial estimate of the column
parameters, we calculate the weights in step~(1a), where
the trimming level $\alpha$ plays the same role as $\alpha^*$ in (I1)
but in a different context. The main purpose of step~1 is to increase
efficiency of the column parameter estimates over those from step~0,
but the corresponding estimates of the row effects might not be robust.
The purpose of step~2 is to robustify the row effect estimates.

General weight functions of $\|\underline{\hat{e}}_i\|^2$ can be
considered in lieu of (\ref{weight}), but we expect that the results
given in the Appendix~\ref{app} still hold under appropriate regularity
conditions. Our proposed robust estimates of parameters $\underline
{\theta}_{k}^{(0)}$ and $\underline{\phi}_k^{(0)}$ are obtained\vadjust{\goodbreak} by
minimizing (\ref{robobj2}) and (\ref{robobj1}), respectively. By
considering the regular SVD on the approximation matrix $\sum_{k=1}^r
\underline{\tilde{\theta}}_k\underline{\tilde{\phi}}{}_k^T$, we actually
obtain a robust SVD on the data matrix $\mathbf{Y}$.

\section{Asymptotic properties}\label{estimates}
The data matrix $\mathbf{Y}$ often arises with the rows representing
individuals randomly sampled from a
large population, but the columns for measurements at $m$ different
locations or time points. It is then
natural to use
$\underline{\theta}$ as the random row effects, and $\underline{\phi}$
as the fixed column effects.
Individuals can be characterized by the row effects, and their spatial
or temporal profiles can be understood by
the column effects. The distinction between the random and the fixed
effects is not relevant
to the optimization problems (\ref{robobj1}) and (\ref{robobj2})
themselves, but is important for the statistical properties of the
estimates obtained from the optimization. To derive the statistical
representations of the row- and column-effect estimates, we use
conditions (M1)--(M5) detailed in Appendix~\ref{AssMod}. Those
conditions also ensure proper parameter identifications.

Following Definition~1.1 of \citet{feng}, we use the rank of the mean
matrix $E(\mathbf{Y})$ as the dimensionality of the model. A
unidimensional model refers to the mean matrix of rank-one. For
unidimensional data, we can use the first
singular component to summarize the row and column effects. For
example, if a unidimensional test of
$m$ items is given to $n$ examinees, the data matrix as the scores of
the examinees on each of the items might be
expected to be of rank one, where a rank-one approximation uses $\theta
_i$ to summarize the ``ability'' of the $i$th examinee
and $\phi_j$ to represent the difficulty level of the $j$th item. In
educational measurements, different forms of unidimensionality has
been used. For a related article on assessing unidimensionality of
polytomous data, see \citet{nandakumar}.

\subsection{Profiling in optimization and column effect estimates}
The number of the~$\theta$'s involved in the objective function (\ref
{robobj1}) increases with $n$, which inconveniences the asymptotic
analysis as $n \to\infty$.\vadjust{\goodbreak}
To bypass this difficulty, we view $\underline{\theta}_k$ as nuisance
parameters in the following profiling procedure. First, we minimize
the objective function (\ref{robobj1}) with respect to
$\underline{\theta}_k $ as
if $\underline{\phi}_k$ ($k=1, \ldots, r$) were known. Then, with
the estimates
$\theta_{ki}^*=\underline{\phi}_k^T\underline{y}_i$, minimizing
(\ref
{robobj1}) is equivalent to minimizing the following objective function:
%
\begin{equation}
\label{robobj3} \min_{\underline{\phi}}\sum_{i=1}^n
\hat{w}_i \Biggl\llVert \Biggl(I_m-\sum
_{k=1}^r\underline{\phi}_k\underline{
\phi}_k^T \Biggr)\underline {y}_i \Biggr
\rrVert ^2,
\end{equation}
under the restrictions that $\|\underline{\phi}_1\|=\cdots
=\|\underline
{\phi}_r\|=1$, and $\underline{\phi}_k\,\bot\,\underline{\phi}_l$ for
$k\neq l$.

Let $\varphi_0=(\underline{\phi}_1^{(0)T},\ldots,\underline{\phi}_r^{(0)T})^T$, ${\vartheta}_{0}=(0,\varphi_0^T)^T$, and $\hat
{\vartheta
}_\tau=(\hat{\xi}_\tau-\xi_\tau,\hat{\varphi}^T)^T$, where $\xi
_\tau$
is the $\tau$th quantile of $\Vert\underline{e}_i \Vert^2$,
$\underline
{e}_i=(I_m-\sum_{k=1}^r\underline{\phi}_k^{(0)}\underline{\phi
}_k^{(0)T})\underline{y}_i$
and $\hat{\varphi}$ is the initial estimate of $\varphi_0$. We obtain
the Bahadur representation for the estimates $\tilde{\varphi
}=(\underline{\tilde{\phi}}{}_1^{T},\ldots,\underline{\tilde{\phi
}}{}_r^{T})^T$ from step~1.
%
\begin{theorem}\label{AsympNorm}
Assume model (\ref{equ0}) with $\hat{\varphi}$ as any root-$n$
consistent estimate of the parameter vector $\varphi_0$. If conditions
\textup{(M1)--(M5)} and \textup{(E1)--(E3)} in Appendix~\ref{ass} hold, then
%
\begin{eqnarray}\label{baha}
\tilde{\varphi}-\varphi_0 &=&-(n\mathbf{D}_0)^{-1}\sum
_{i=1}^n w_i \pmatrix{
\underline{b}_1 \bigl(\theta_{1i}^{(0)},\ldots,
\theta _{ri}^{(0)},\underline{\varepsilon}_i,
\varphi_0 \bigr)
\vspace*{3pt}\cr
\vdots
\vspace*{3pt}\cr
\underline{b}_r \bigl(
\theta_{1i}^{(0)},\ldots,\theta _{ri}^{(0)},
\underline{\varepsilon}_i,\varphi_0 \bigr)}
\nonumber\\[-8pt]\\[-8pt]
&&{}+\mathbf{G}_n^T
\pmatrix{ \hat\vartheta_{1-\alpha}-\vartheta_{0}
\vspace*{3pt}\cr
\hat
\vartheta_{\alpha}-\vartheta_{0}} +o_p
\bigl(n^{-1/2} \bigr),\nonumber
\end{eqnarray}
where $w_i = 1(\xi_\alpha<\|\underline{e}_i\|^2\leq\xi_{1-\alpha
}) $,
$\mathbf{D}_0$ is an $mr \times mr$ nonsingular square matrix,
$\mathbf
{G}_n$ is an $mr\times2(mr+1)$ matrix with the Frobenius norm $\|
\mathbf{G}_n\|_F=O(1)$ and
\begin{eqnarray*}
&& \underline{b}_j \bigl(\theta_{1i}^{(0)},\ldots,
\theta _{ri}^{(0)},\underline {\varepsilon}_i,
\varphi_0 \bigr)
\\
&&\qquad =2 \bigl\{\theta_{ji}^{(0)}+
\underline{\varepsilon}_i^T\underline{\phi
}_j^{(0)} \bigr\}^2\underline{
\phi}_j^{(0)}- \bigl\{\theta_{ji}^{(0)}+
\underline {\varepsilon}_i^T\underline{
\phi}_j^{(0)} \bigr\}\underline{y}_i -\sum
_{k=1}^r \bigl\{\theta_{ki}^{(0)}+
\underline{\varepsilon}_i^T\underline {
\phi}_k^{(0)} \bigr\}\underline{y}_i.
\end{eqnarray*}
\end{theorem}

The specific forms of $\mathbf{D}_0$ and $\mathbf{G}_n$ can be found in
the supplementary material [\citet{feng1}]. From Theorem~\ref{AsympNorm},
Lemma~\ref{quantile} (in the Appendix) and Theorem~2.2 of \citet{feng},
it is clear that the estimate $\tilde{\varphi}$ of the parameter vector
$\varphi_0$ is root-$n$ consistent with asymptotic normality. Its
asymptotic variance--covariance matrix is complicated because both
variations from the initial estimates and the variation from the
weighted least squares method are present.

\subsection{Row effect predictions}
Note that the least squares estimate of $\theta_{ki}^{(0)}$ is~$\underline{\tilde\phi}{}_k^T\underline{y}_i$, so it can be seriously
affected\vadjust{\goodbreak} by any outlying value of the observed vector $\underline
{y}_i$. We now consider the robust procedure that minimizes (\ref
{robobj2}) for a smooth loss function $L$.

If $L$ has continuous second derivative, the minimizers of (\ref
{robobj2}) are, by the implicit function theorem in calculus,
%
\begin{eqnarray}
\label{implicit1} {\tilde\theta}_{1i}&=&f(\underline{y}_i,
\underline{\tilde\phi}_1,\underline{\tilde\phi}_2,
\ldots,\underline{\tilde\phi}_r),
\\
\nonumber &&{} \vdots
\\
\label{implicit2} {\tilde\theta}_{ri}&=&f(\underline{y}_i,
\underline{\tilde\phi}_r,\underline{\tilde\phi}_1,
\ldots,\underline{\tilde\phi}_{r-1}),
\end{eqnarray}
where $f$ is a function with continuous partial derivatives with
respect to $\phi_{kj}$ for $k=1,\ldots,r$ and $j=1,\ldots,m$.

Before we move on, it helps to explore some properties of the implicit
function~$f$.
Consider minimizing the following objective function:
\[
\sum
_{j=1}^mL \Biggl(y_{ij}-\sum
_{k=1}^r\theta_{ki}{\phi
}^{(0)}_{kj} \Biggr),
\]
which can be written under model (\ref{equ0}) as
\[
\sum
_{j=1}^mL \Biggl(\varepsilon_{ij}-
\sum_{k=1}^r \bigl(\theta _{ki}-
\theta^{(0)}_{ki} \bigr)\phi^{(0)}_{kj}
\Biggr).
\]
When this minimization is performed with respect to $\theta_{ki}$, we have
%
\begin{eqnarray}
\label{relation1} f \bigl(\underline{y}_i,\underline{
\phi}{}^{(0)}_1,\underline{\phi}{}^{(0)}_2,
\ldots,\underline{\phi}{}^{(0)}_r \bigr) &=&
\theta^{(0)}_{1i}+f \bigl(\underline{\varepsilon}_i,
\underline{\phi}{}^{(0)}_1, \underline{\phi}{}^{(0)}_2,
\ldots,\underline{\phi}{}^{(0)}_r \bigr),
\\
\nonumber &&{} \vdots
\\
\label{relation2} f \bigl(\underline{y}_i,\underline{
\phi}{}^{(0)}_r,\underline{\phi}{}^{(0)}_1,
\ldots,\underline{\phi}{}^{(0)}_{r-1} \bigr) &=&
\theta^{(0)}_{ri}+f \bigl(\underline{\varepsilon}_i,
\underline{\phi}{}^{(0)}_r,\underline{\phi}{}^{(0)}_1,
\ldots,\underline{\phi}{}^{(0)}_{r-1} \bigr).
\end{eqnarray}
If $L$ is even, then the function $f$ is radially symmetrical with
respect to its first argument.
We obtain the asymptotic result for the estimates $\tilde\theta_{ki}$
defined as the minimizer of (\ref{robobj2}) in the following theorem.
%
\begin{theorem}\label{Consistrow}
Assume model (\ref{equ0}) with $\hat{\varphi}$ as any root-$n$
consistent estimate of the parameter vector $\varphi_0$. If conditions
\textup{(M1)--(M5)}, \textup{(A1)--(A4)} and \textup{(C3)} in Appendix~\ref{ass} hold, then
\begin{eqnarray*}
\sum_{k=1}^r\tilde{\theta}_{ki}
\underline{\tilde{\phi}}_k
&\stackrel {d} {\longrightarrow}&\sum
_{k=1}^r\theta_{ki}^{(0)}
\underline{{\phi}}_k^{(0)}
\\
&&{}+f \bigl(\underline{\varepsilon}_i,\underline{
\phi}{}^{(0)}_1,\underline {\phi}{}^{(0)}_2,
\ldots,\underline{\phi}{}^{(0)}_r \bigr)\underline{{\phi
}}_1^{(0)}
\\
&&{} +\cdots+f \bigl(\underline{\varepsilon}_i,
\underline{\phi}{}^{(0)}_r,\underline{\phi}{}^{(0)}_1,
\ldots,\underline{\phi}{}^{(0)}_{r-1} \bigr)\underline{{
\phi}}_r^{(0)},
\end{eqnarray*}
where $\tilde{\theta}_{ki}$ is defined in (\ref{implicit1})--(\ref
{implicit2}), and $\stackrel{d}{\longrightarrow}$ refers to convergence
in distribution.
\end{theorem}

It is clear from Theorem~\ref{Consistrow} that each row of the
approximating matrix $\sum_{k=1}^r \underline{\tilde{\theta
}}_k\underline{\tilde{\phi}}{}_k^T$ converges\vspace*{1pt} in distribution to the
corresponding row of the rank-$r$ matrix $\sum_{k=1}^r\underline
{\theta
}{}^{(0)}_k\underline{\phi}_k^{(0)T}$ and some function of the model
errors $\underline{\varepsilon}$.

\section{Application}\label{appli}
For vector measurements, a unidimensional summary is widely used in
data analysis. In this section, we consider testing on the sufficiency
of unidimensional summaries, against the alternative that the matrix
$\mathbf{Y}$ is a rank two matrix under model (\ref{equ0}).

\subsection{Hypothesis testing}
With the asymptotic results of the previous section, we consider
hypothesis testing here based on the robust estimates.
The null hypothesis is $\underline{\mu}_2=\underline{0}$, which implies
unidimensionality of the mean matrix $E(\mathbf{Y})$, and that no
meaningful pattern can be found in the second dimension of the data matrix.
This hypothesis is especially interesting in the probe-level microarray
data analysis, where unidimensional models are usually assumed to
summarize the gene expression level from the intensity data matrix
[\citet{LiAndWong}, \citet{irizarry}].\vadjust{\goodbreak}

We first consider the estimation by minimizing (\ref{robobj1}) with
$r=2$. We then use the column vectors $\tilde{\underline{\phi}}_1$ and
$\underline{0}$ in minimizing (\ref{robobj2}) to obtain the estimate
$\tilde{\theta}_{1i}=f(\underline{y}_i,\tilde{\underline{\phi
}}_1,\underline{0})$, where $f$ is defined in (\ref
{implicit1})--(\ref
{implicit2}). For convenience, we use $f(\underline{y}_i,\tilde
{\underline{\phi}}_1)$ instead of $f(\underline{y}_i,\tilde
{\underline
{\phi}}_1,\underline{0})$ from now on.
Let
\[
\gamma(\underline{y}_i,\varphi)=\sum_{j=1}^mL'
\bigl(y_{ij}-f(\underline {y}_i,{\underline{
\phi}_1}){\phi}_{1j} \bigr){\phi}_{2j}
\]
be the score for unidimensionality corresponding to the $i$th vector
$\underline{y}_i$.
We have the following result.
%
\begin{theorem}\label{sufficiency}
Let $\underline{a}=(a_1,\ldots,a_n)^T$
be a vector that is orthogonal to $\underline{\mu}_1$ and satisfies
$\|\underline{a}\|^2=n$ with a bounded supremum norm. Assume model
(\ref{equ0}) and conditions \textup{(M1)--(M5)}, \textup{(C1)--(C4)}, \textup{(D1)--(D2)} in
Appendix~\ref{ass}, then
%
\begin{equation}
\label{asym1} {n^{-1/2}\underline{a}{}^T\underline{\tilde{
\gamma}}}/\tilde {\sigma}_n\stackrel{L} {\longrightarrow}N(0,1),
\end{equation}
under the null hypothesis that $\underline{\mu}_2=\underline{0}$, where
%
\begin{eqnarray}
\label{gamma} \underline{\tilde{\gamma}} &=& \bigl(\gamma(\underline{y}_1,
\tilde \varphi ),\ldots,\gamma(\underline{y}_n,\tilde\varphi)
\bigr)^T,
\\
\label{robsigma} \tilde{\sigma}_n^2 &=& n^{-1}\sum
_{i=1}^n\gamma^2(\underline
{y}_i,\tilde \varphi)- \Biggl\{n^{-1}\sum
_{i=1}^n\gamma(\underline{y}_i,\tilde
\varphi ) \Biggr\}^2
\end{eqnarray}
and $\tilde{\varphi}$ is the robust estimate defined in Section~\ref{init}.
\end{theorem}
%
%
\begin{remark}
If the loss function $L$ is the $L_2$ norm, then $L'(x)=2x$. It then
follows that\vspace*{1pt}
$\gamma(\underline{y}_i,\tilde\varphi)
=2\sum_{j=1}^m\{y_{ij}-(\tilde{\underline{\phi}}{}_1^T\underline
{y}_i)\tilde{\phi}_{1j}\}\tilde{\phi}_{2j}
=2\tilde{\underline{\phi}}{}_2^T\underline{y}_i$,\vadjust{\goodbreak} because $\tilde
{\underline{\phi}}_1\,\bot\,\tilde{\underline{\phi}}_2$ for the least
squares case.
Thus, the statistic used by \citet{feng} can be viewed as a special case
of Theorem~\ref{sufficiency}.
\end{remark}

If the direction vector $\underline{a}$ is not orthogonal to
$\underline
{\mu}_1$, then $n^{-1/2}\underline{a}{}^T\underline{\tilde{\gamma
}}
/\tilde{\sigma}_n$ may not converge in distribution to a mean
zero distribution.
Typically $\underline{\mu}_1$ is unknown and needs to be estimated.
This is usually done by extra group information in the rows to enable
us to consistently estimate $\underline{\mu}_1$, which is sufficient to
have the asymptotic result for the pivotal statistic in Theorem~\ref
{sufficiency}. This theorem also ensures the validity of the bootstrap
as described in Section~3.3 of \citet{feng} based on Theorem~1 of \citet{Mammen}.

It is certainly possible that the direction vector $\underline{a}$
happens to be a poor choice in the sense of low power against a
particular alternative. To ensure decent power of the test, we can
consider several target directions that are orthogonal to each
other.
%
\begin{theorem}\label{test2}
Assume the conditions of Theorem~\ref{sufficiency}.
Consider a $K\times n$ matrix $A$ with all the row vectors
orthogonal to each other, with $K$ being fixed. If the vector
$\underline{a}_l=(a_{l1},\ldots,a_{ln})^T$ is the $l$th row of the
matrix $A$ and satisfies $\underline{a}_l\,\bot\,\underline{\mu}_1$,\vadjust{\goodbreak} and
$\|\underline{a}_l\|^2=n$ with uniformly bounded elements,
then
$P(n^{-1}\|A\underline{\tilde\gamma}\|^2/\tilde{\sigma}^2\leq
x)-F_K (x)\rightarrow0$ under the null hypothesis that
$\underline{\mu}_2=0$, where $F_K$ is the cumulative distribution
function of the $\chi^2_K$ distribution, $\tilde{\underline{\gamma}}$
and $\tilde{\sigma}$ are given
in (\ref{gamma}) and (\ref{robsigma}), respectively.
\end{theorem}

\subsection{A simulation study}\label{simulation}
In this section, we use a simulation study to assess the performance of
the target direction test based on robust loss functions. We
independently generate 20 rows of size 12 from model (\ref{equ0}), with
the mean of the corresponding $20\times12$ matrix equal to $\underline
{\mu}_1\underline{\phi}_1^T$ and $\underline{\mu}_1\underline
{\phi
}_1^T+\underline{\mu}_2\underline{\phi}_2^T$ under the null and the
alternative hypotheses, respectively,
where $\underline{\mu}_1=(20,\ldots,20)^T$, $\underline{\mu
}_2=2^{1/2}(1,-1,\ldots,1,-1)^T$, $\underline{\phi}_1=(1,\ldots,1)^T/12^{1/2}$ and $\underline{\phi}_2=(1,-1,\ldots,1,-1)^T/12^{1/2}$. The random effects $\theta_{1i}^{(0)}-\mu_{1i}$ and
$\theta_{2i}^{(0)}-\mu_{2i}$ are generated from normal distributions
with mean 0 and variances 4 and 1, respectively.

To assess the robustness of the method, we generate model errors in two
ways. In an outlier-free model, all the errors are independently
generated from one of the three cases: (I) $2^{-1/2} N(0,1)$; (II)
$(3/10)^{-1/2}t_5$, where $t_5$ is the $t$ distribution with 5 degrees
of freedom; (III) $2^{-1}(\chi^2_1-1)$, where $\chi^2_1$ is the $\chi
^2$ distribution with 1~degree of freedom. In a contaminated model, the
first two rows of the matrix are generated from the mixture of the
normal distribution $N(0,11)$ with probability 0.1 and one of the three
distributions (I), (II) or (III) with probability 0.9, but the other
rows are generated as in the outlier-free model. Under the contaminated
model, outliers are likely to occur in the first two rows. A total of
5000 data sets are generated from each model in the simulation
study.\vadjust{\goodbreak}

For the initial steps (I1)--(I5) of Section~\ref{init}, we use $\alpha
^*=0.3$ and 100 randomly selected subsets, and the constant $\alpha
=0.1$ is used in calculating the weights (\ref{weight}). With only two
possible outlying rows, the probability that all 100 subsets contain an
outlier is less than 0.001.

We consider two choices of the direction vector $\underline{a}$, with
$\underline{a} \propto\underline{\mu}_2$ in the first case, and
$\underline{a} \propto(3/2)^{1/2}(1,-1,\ldots,1,-1)^T+(1,\ldots,1,-1,\ldots,-1)^T$ in the \mbox{second}~case. The bootstrap
calibration method of \citet{feng} is used to calculate the $p$-values
of the tests. Three loss functions are used for comparison. They~are
\begin{longlist}[(L3)]
\item[(L1)]``Logistic'':
$L(s)=C\log (\cosh(s/C) )$,
\item[(L2)] ``Huber'':
\[
L(s)=\cases{2^{-1}s^2,&\quad $|s|\leq C$,
\vspace*{3pt}\cr
C|s|-2^{-1}C^2,&\quad $|s|>C$,}
\]
\item[(L3)]``Least squares'': $L(s)=s^2$,
\end{longlist}
where $C=0.1$ is used in our simulation. Since $C$ is close to zero,
the two robust loss functions (L1) and (L2) lead to results that are
similar to those obtained under the $L_1$ loss $L(s)=|s|$.

\begin{table}[t]
\tabcolsep=0pt
\caption{Estimated type~I errors and powers of various tests at the
nominal level of 5\%, with data generated from outlier-free
models}\label{sim2}
\begin{tabular*}{\tablewidth}{@{\extracolsep{\fill}}lcccccc@{}}
\hline
&\multicolumn{3}{c}{\textbf{Null}}&\multicolumn{3}{c@{}}{\textbf{Alternative}}\\[-6pt]
&\multicolumn{3}{c}{\hrulefill}&\multicolumn{3}{c@{}}{\hrulefill}
\\
\textbf{Size}&\textbf{Normal} &$\bolds{t}$ & $\bolds{\chi^2}$& \textbf{Normal}&$\bolds{t}$&$\bolds{\chi^2}$
\\
\hline
Logistic\tabnoteref[a]{t1a}&0.051&0.049&0.043&1.000&0.999&0.995 \\
Huber\tabnoteref[a]{t1a}&0.051&0.049&0.043&1.000&0.999&0.997\\
Least squares\tabnoteref[a]{t1a}&0.050&0.045&0.035&1.000&0.999&0.998 \\
Logistic\tabnoteref[b]{t1b}&0.052&0.054&0.051&0.941& 0.959&0.978 \\
Huber\tabnoteref[b]{t1b}&0.053&0.054&0.051&0.936&0.956&0.976\\
Least squares\tabnoteref[b]{t1b}&0.054&0.046&0.039&1.000&0.989&0.998 \\
\hline
\end{tabular*}
\tabnotetext[a]{t1a}{The results are from the case where $\underline{a}\propto\underline{\mu}_2$.}
\tabnotetext[b]{t1b}{The results are from the case where $\underline{a}\propto(3/2)^{1/2}(1,-1,\ldots,1,-1)^T+(1,\ldots,1,\break -1,\ldots,-1)^T$.}
\end{table}

\begin{table}
\caption{Estimated type~I errors and powers of various tests at the
nominal level of 5\%, with data generated from contaminated models}\label{sim1}
\begin{tabular*}{\tablewidth}{@{\extracolsep{\fill}}lcccccc@{}}
\hline
&\multicolumn{3}{c}{\textbf{Null}}&\multicolumn{3}{c@{}}{\textbf{Alternative}}\\[-6pt]
&\multicolumn{3}{c}{\hrulefill}&\multicolumn{3}{c@{}}{\hrulefill}
\\
\textbf{Size}&\textbf{Normal} &$\bolds{t}$ & $\bolds{\chi^2}$& \textbf{Normal}&$\bolds{t}$&$\bolds{\chi^2}$
\\
\hline
Logistic\tabnoteref[a]{t2a}&0.049&0.051&0.052&0.987&0.983&0.985 \\
Huber\tabnoteref[a]{t2a}&0.049&0.048&0.052&0.987&0.983&0.986\\
Least squares\tabnoteref[a]{t2a}&0.024&0.021&0.019&0.467&0.398&0.404 \\
Logistic\tabnoteref[b]{t2b}&0.054&0.046&0.053&0.884& 0.908&0.866 \\
Huber\tabnoteref[b]{t2b}&0.054&0.048&0.052&0.882&0.906&0.862\\
Least squares\tabnoteref[b]{t2b}&0.021&0.018&0.022&0.455&0.371&0.464\\
\hline
\end{tabular*}
\tabnotetext[a]{t2a}{The results are from the case where $\underline{a}\propto\underline{\mu}_2$.}
\tabnotetext[b]{t2b}{The results are from the case where $\underline{a}\propto(3/2)^{1/2}(1,-1,\ldots,1,-1)^T+(1,\ldots,1,\break -1,\ldots,-1)^T$.}
\vspace*{4pt}
\end{table}

We summarize the results for the outlier-free models in Table~\ref{sim2}.
It is clear from Table~\ref{sim2} that all the three tests preserve
type~I errors well, and they achieve very high power under the
alternative. The story is different, however, for the contaminated
models with the results in Table~\ref{sim1}. When no more than 10\% of
outliers are present, the test based on the square loss becomes too
conservative with low power, but the robust tests with (L1) and (L2)
loss functions withstand the outliers very well.

\subsection{Case study}\label{application}
In this section, we analyze a real microarray dataset and
examine the test results based\vadjust{\goodbreak} on the least squares method of \citet{feng} as well as the robust alternative
studied in this paper. We use the same GeneChip
data obtained from the MicroArray Quality Control project [\citet
{shi}, \citet{Guixian}].
There are a total of 20 microarrays
(HG-U133-Plus-2.0) with 54,675 probe-sets (each composed of 11 probes)
on each, generated from
five colorectal adenocarcinomas and five matched normal colonic
tissues with one technical replicate at each of two laboratories
involved in the MAQC project. We use the intensity measure of perfect matches,
and preprocess the probe-level microarray data with the
``RMA'' background adjustment method [\citet{irizarry}] and the quantile
normalization
method [\citet{bolstad}].


We consider a target direction [see supplementary material, \citet
{feng1}] to contrast the
two groups: the normal tissue group and the tumor group. Since the gene
expressions from the arrays of the same group
are expected to be equal, the target direction is approximately
orthogonal to the mean of ${\underline\theta}_1$.

%
\begin{figure}

\includegraphics{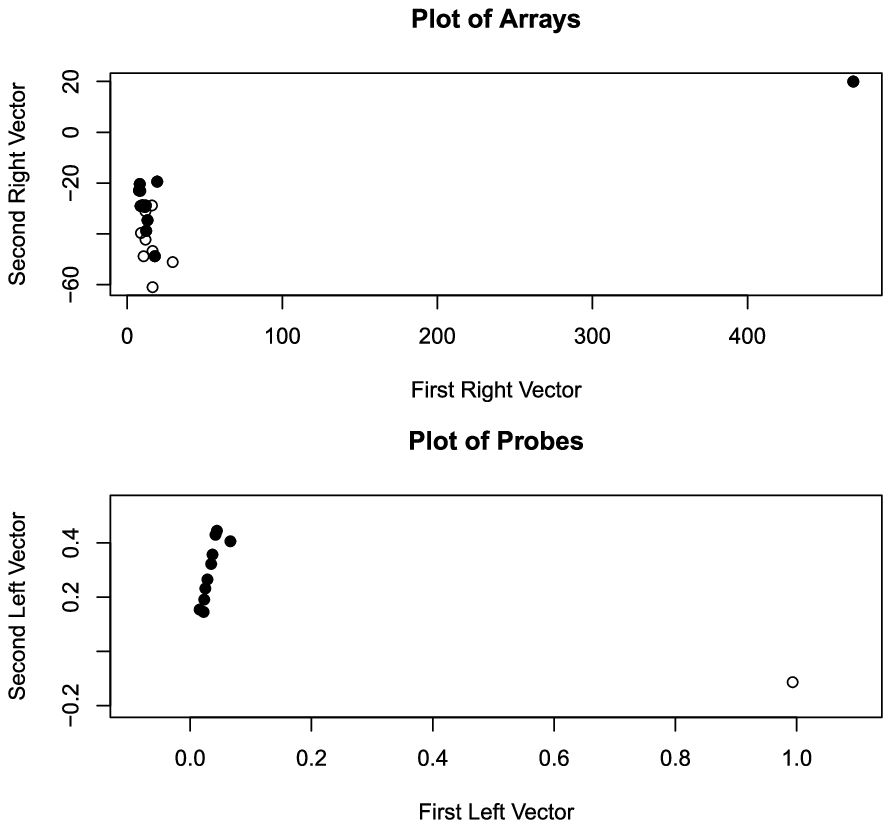}

\caption{Scatter plot of singular vectors for the
probe-set ``1555106\_a\_at'' from the regular SVD. In the upper panel,
the circles represent the arrays from tumor samples, while the~solid
points represent normal tissues. In the lower panel, the circle
corresponds to probe~3.\looseness=1}\label{1555106}\vspace*{4pt}
\end{figure}

For the first approximation by minimizing (\ref{robobj1}), we use the
same values of $\alpha$~and~$\alpha^*$ as those of Section~\ref{simulation}.
For the second approximation by minimizing~(\ref{robobj2}),
we consider two loss functions: one is for the square loss and the
other is
the logistic loss function with $C=1.205$ (times the scale of the
residuals). With this choice of $C$, we retain $95\%$ asymptotic
efficiency at the normal distribution.




We inspect one probe-set ``1555106\_a\_at'' to better understand the
discrepancies between the least squares method and the robust
alternative. In this case, the data matrix has 20 rows and $m=11$
columns. In Figure~\ref{1555106}, we plot the arrays and the probes
with the coordinates $(\hat{\theta}_{1i},\hat{\theta}_{2i})$ and
$(\hat
{\phi}_{1j},\hat{\phi}_{2j})$, respectively, for $i=1,\ldots,20$ and
$j=1,\ldots,11$, where the least squares estimates are used.
The \mbox{$p$-}value is $0.036$ based on the
least squares method, and the first four singular values are $(472,
163, 36, 29)$. It is clear from Figure~\ref{1555106} that there exist
an outlying array and an outlying probe.
Further inspection of the data shows that there exists
an outlying measurement in the
outlying array and the outlying probe in the intensity data matrix.
In other words, it is likely that the significant two-dimensional mean
structure is caused by the outlier.

\begin{figure}

\includegraphics{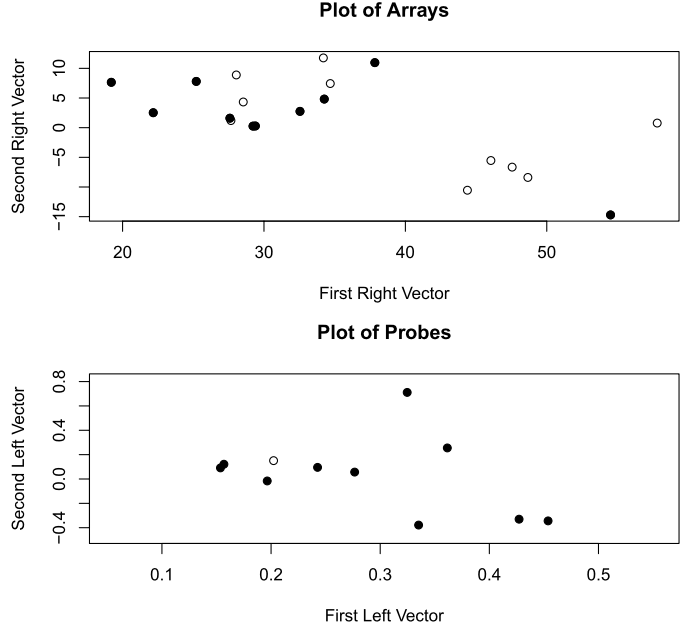}

\caption{Scatter plot of singular vectors for the
probe-set ``1555106\_a\_at'' from a robust approximation. In the upper
panel, the circles represent the arrays from tumor samples, while the
solid points represent normal tissues. In the lower panel, the circle
corresponds to probe~3.}\label{15551061}
\end{figure}

With the robust alternative, the $p$-value is $0.741$, and no outlying
estimates of the arrays-effects or probe-effects are
observed in Figure~\ref{15551061}.
The first four singular
values are $(169,29,25,23)$ in this case, and the second singular value
is close to the third and the fourth, which
indicates that the 2nd singular structure is likely to be due to noise.
From this empirical example, we see that the robust method is
powerful in moderating the effect from outliers.
More details of the case study can be found in the supplementary
material [\citet{feng1}].\newpage

%



\begin{appendix}
\section{Assumptions}\label{ass}

\subsection{Model assumptions}\label{AssMod}
\begin{longlist}[(M1)]
\item[(M1)] The column vectors $\underline{\phi}_k^{(0)}$ $(k=1,
\ldots, r$) are orthogonal to one another.

\item[(M2)] The row vectors $\underline{\theta}_{k}^{(0)}$ $(k=1,
\ldots, r$) are
independently distributed with mean
$\underline{\mu}_k=(\mu_{k1},\ldots,\mu_{kn})^T$ and variance
$\sigma_k^2I_n$, for $k=1,\ldots,r$. The components of $\underline
{\theta}_k^{(0)}$
are independently distributed with finite fourth moments.
Moreover, $\underline{\mu}_k\,\bot\,\underline{\mu}_l$, for $k\neq l$,
where $\bot$ denotes orthogonality.
\item[(M3)]The error variables $\underline{\varepsilon
}_i=(\varepsilon
_{i1},\ldots,\varepsilon_{im})^T$
are independently generated from a distribution with mean zero and
finite fourth moment, and $\varepsilon_{ij}$ is symmetrically
distributed with $E(\varepsilon_{ij}^2)=\sigma^2$.
%
\item[(M4)] The variables $\{\theta_{1i}^{(0)}\},\ldots,\{\theta_{ri}^{(0)}\}$ and $\{\underline{\varepsilon}_i\}$ are
mutually independent.
\item[(M5)]$n^{-1}\|\underline{\mu}_{k}\|^2\rightarrow\mu_k^2$
as $n\rightarrow\infty$ for some finite constants $\mu_k$, where
$\|\cdot\|^2$ is the $L_2$ norm. We assume that
$\mu_k^2+\sigma^2_k>\mu_l^2+\sigma_l^2$ when $k<l$, which is necessary
for the
identifiability of the model parameters.
\end{longlist}
%

These assumptions are clearly satisfied with Gaussian row-effects and
Gaussian errors. In robust statistics, a traditional parametric model
is often assumed for the outlier-free part of the data, but we design a
robust procedure to be insensitive to data contamination.

\subsection{\texorpdfstring{Assumptions for Lemma \protect\ref{quantile} and Theorem \protect\ref{AsympNorm}}
{Assumptions for Lemma B.2 and Theorem 3.1}}\label{Asymp1}

Let $\vartheta=(\delta,\varphi^T)^T$, and $\hat{\xi}_{\tau}$ and
${\xi
}_{\tau}$ be the sample and the\vspace*{1pt} population $\tau$ quantiles of
$\|\underline{\hat{e}}_i\|^2$ and $\|\underline{e}_i\|^2$,
respectively, where $\underline{\hat{e}}_i= (I_m-\sum_{k=1}^r\underline{\hat{\phi}}_k\underline{\hat{\phi}}{}_k^T
)\underline{y}_i$ and $\underline{e}_i= (I_m-\break \sum_{k=1}^r\underline
{{\phi}}_k^{(0)}\underline{{\phi}}_k^{(0)T} )\underline{y}_i$.
Let the function $g$ denote the probability density of the random
variable $\|\underline{e}_i\|^2$.
\begin{longlist}[(E1)]
%
\item[(E1)]
The value $g(\xi_\tau)$ is bounded and positive, and $g$ is continuous
in a neighbour of $\xi_\tau$.
%
\item[(E2)] 
$n^{-1}\sum_{i=1}^nE \{\|\underline{y}_i\|\llvert \frac{\partial
{f_i}}{\partial\underline{y}_i}\rrvert  /f_i(\underline
{y}_i) \}
\leq K$ for some constant $K$ and all $n$, where $f_i$ is the
probability density function of the random vector $\underline{y}_i$.
\item[(E3)] For given $\xi\in\mathbb{R}$, $n^{-1}\sum_{i=1}^n\mathbf
{H}_i(\xi,\underline{\phi}_j^{(0)},\vartheta_0)=O(1)$, for
$j=1,\ldots,r$, where
%
\begin{equation}
\label{Hi} \mathbf{H}_i(\xi,\underline{\nu},\vartheta)=
\frac{\partial E\{
\mathbf
{M}(\xi,\delta,\varphi,\underline{y}_i)\underline{\nu}\}
}{\partial
\vartheta}
\end{equation}
and
%
\begin{equation}
\label{Mn} \mathbf{M}(\xi,\delta,\varphi,\underline{y}_i)=1
\Biggl\{ \Biggl\llVert \Biggl(I_m-\sum_{k=1}^r
\underline{{\phi}}_k\underline{{\phi}}_k^{T}
\Biggr)\underline{y}_i \Biggr\rrVert ^2\leq{\xi}+\delta
\Biggr\}\underline {y}_i\underline{y}_i^T.
\end{equation}
\end{longlist}

%
\begin{remark}
By similar arguments to those used in the proof of Lem\-ma~\ref
{quantile}, we note that assumption (E3) holds if $n^{-1}\sum_{i=1}^nE \{\|\underline{y}_i\|^3\llvert \frac{\partial
{f_i}}{\partial\underline{y}_i}\rrvert/\break f_i(\underline
{y}_i) \}
\leq K'$ for some constant $K'$ and all $n$. Conditions (E2) and (E3)
are satisfied by the Gaussian distribution as well as any t
distribution with finite fourth moment.
\end{remark}

\subsection{Assumptions on the loss function}\label{AssLemma1}
\begin{longlist}[(C1)]
\item[(C1)]The loss function $L$ is even and nonnegative, and $L(x)=0$
if and only if $x=0$.
\item[(C2)]The first derivative $L'$ is continuous, piecewise
differentiable, nondecreasing in $\mathbb{R}$ and positive in
$\mathbb{R}^+$.
\item[(C3)]The second derivative $L''$ is nonnegative, nonincreasing in
$\mathbb{R}^+$ and piecewise continuous.
\item[(C4)]The derivatives $L'$ and $L''$ satisfy $|L'(x)|\leq C_0|x|$
and $L''(x)\leq C_0$ at all $x\in\mathbb{R}$, for some constant $C_0$.
\end{longlist}


\subsection{\texorpdfstring{Assumptions for Theorem \protect\ref{sufficiency}}
{Assumptions for Theorem 4.1}}\label{Small}
\begin{longlist}[(D2)]
\item[(D1)]$\max_{1\leq i\leq n}\|\underline{y}_i\|
=O_p(n^{1/4-\delta
})$ for some small positive number $\delta$.
\item[(D2)] The distribution of $\theta_{2i}^{(0)}-\mu_{2i}$ is
symmetric around zero.
\end{longlist}
%

\section{Proofs}\label{app}
In the proofs, we assume that $r=2$ for simplicity. The same arguments
work for the general cases of $r \ge2$.
We first give the Bahadur representations of the quantile estimates.
First we state four lemmas, but their proofs can be found in the
supplementary material.

%
\begin{lemma}\label{initial}
Suppose that assumptions \textup{(M1)--(M5)} hold and \textup{(I2)--(I4)} in step~0 are
repeated a fixed number of times, then the initial estimate $\hat
{\varphi}$ is \mbox{root-$n$} consistent for $\varphi_0$.
\end{lemma}

%
\begin{lemma}\label{quantile}
Suppose that assumptions \textup{(M1)--(M5)} and \textup{(E1)--(E2)} hold, and $\hat
{\varphi}$ is the initial root-$n$ consistent estimate of $\varphi
_0$, then
%
\begin{equation}
\label{quantrep} \hat\xi_\tau-\xi_\tau=- \bigl\{ng(
\xi_\tau) \bigr\}^{-1}\sum_{i=1}^n
\psi _\tau \bigl\{ \|\underline{e}_i\|^2-
\xi_\tau \bigr\}+O_p \bigl(n^{-1/2} \bigr),
\end{equation}
where $\hat\xi_\tau$, $\xi_\tau$, $\underline{e}_i$, $g$ and $v$ are
defined in Appendix~\ref{Asymp1}, and \mbox{$\psi_\tau(u)=\tau-1(u<0)$}.
\end{lemma}

%
\begin{lemma}\label{app0}
If conditions \textup{(M2)}, \textup{(M3)} and \textup{(M5)} hold, and $\hat{\varphi}$ is the
initial root-$n$ consistent estimate of $\varphi_0$, then
%
\begin{eqnarray}\label{app1}
&&n^{-1}\sum_{i=1}^n1
\bigl\{\hat{\xi}_{\alpha}< \|\underline {\hat{e}}_i\|^2
\leq\hat{\xi}_{1-\alpha} \bigr\}\underline {y}_i
\underline{y}_i^T\nonumber
\\
&&\qquad \stackrel{p} {\longrightarrow}(1-2\alpha) \bigl(\mu_1^2+
\sigma _1^2 \bigr)\underline{{\phi}}_1^{(0)}
\underline{{\phi}}_1^{(0)T}
\\
&&\hspace*{41pt}{} +(1-2\alpha ) \bigl(\mu_2^2+\sigma_2^2 \bigr)
\underline{{\phi}}_2^{(0)}\underline{{\phi}}_2^{(0)T}+
\sigma ^2(\alpha)I,\nonumber
\end{eqnarray}
where $\sigma^2(\alpha)=E [1 \{\xi_\alpha\leq\|(I_m-\sum_{k=1}^r\underline{{\phi}}_k^{(0)}\underline{{\phi
}}_k^{(0)T})\underline
{\varepsilon}_i\|^2\leq\xi_{1-\alpha} \}\varepsilon
_{ij}^2 ]$.
\end{lemma}

Let $\tilde{\varphi}$ be the estimate of $\varphi_0$ from step~1. We
now have
%
\begin{lemma}\label{Consistency}
Suppose that the observations $\underline{y}_i,\underline{y}_2,\ldots,\underline{y}_n$ are
drawn from model (\ref{equ0}). If assumptions \textup{(M1)--(M5)} hold, then
$\tilde{\varphi}\stackrel{p}{\rightarrow}\varphi_{0}$.
\end{lemma}



In the following lemma, we obtain upper bounds for the estimates of
$\underline{\theta}{}^{(0)}$'s given~$\varphi$.
%
\begin{lemma}\label{moment} If conditions \textup{(C1)} and \textup{(C2)} hold, then we have
%
\begin{equation}
\label{upper} f^2(\underline{y},\underline{\phi}_1,
\underline{\phi}_2,\ldots,\underline{\phi}_r)+
\cdots+f^2(\underline{y},\underline{\phi}_r,\underline{
\phi}_1,\ldots,\underline{\phi}_{r-1})\leq
4m^2\|\underline{y}\|^2
\end{equation}
for any $\varphi\in\mathbb{S}$ and $\underline{y} \in\mathbb
{R}^m$ where
%
\begin{equation}
\label{S} \mathbb{S}= \bigl\{ \bigl(\underline{\phi}_1^T,
\ldots,\underline{\phi}_r^T \bigr)^T\in
\mathbb{R}^{rm}\dvtx  \|\underline{\phi}_k\|=1,\underline{
\phi}_k\,\bot\,\underline{\phi}_l,\mbox{ for }k\neq l \bigr\}
\end{equation}
and $f$ is defined in (\ref{implicit1})--(\ref{implicit2}).
Furthermore,
\[
f^2 \bigl(\underline{y},\underline{\phi}^*_1,
\underline{ \phi}^*_2,\ldots,\underline{\phi}^*_r \bigr)+
\cdots+f^2 \bigl(\underline{y},\underline{\phi}^*_r,
\underline{\phi}^*_1,\ldots,\underline{\phi}^*_{r-1} \bigr)
\leq 4m^3(1-3\tau/2)^{-1}\|\underline{y}\|^2,
\]
where $\varphi^*=\lambda\varphi_1+(1-\lambda)\varphi_2$, $\lambda
\in
(0,1)$, and $\|\varphi_1-\varphi_2\|\leq\tau$ for $0<\tau<2/3$ and
$\varphi_1,\varphi_2\in\mathbb{S}$.
\end{lemma}
%
%
\begin{remark}
The result of Lemma~\ref{moment} holds uniformly for $\varphi\in
\mathbb
{S}$, so the existence of the moments of $\underline{y}_i$ ensures the
existence of the corresponding moments of the estimates of $\underline
{\theta}{}^{(0)}$'s given $\varphi$.
\end{remark}

\begin{pf*}{Proof of Lemma~\ref{moment}}
Again we present the proof for $r=2$.
From the definition of $f$, we have
\[
\sum_{j=1}^mL \bigl(y_{j}-f(
\underline{y},\underline{\phi}_{1},\underline {\phi}_2){
\phi}_{1j}-f(\underline{y},\underline{\phi}_{2},\underline {
\phi}_1)\phi_{2j} \bigr)\leq\sum
_{j=1}^mL(y_{j}),
\]
where $y_j$ is the $j$th component of any vector $\underline{y} \in R^m$.

From condition (C1), we have
\[
L \bigl(y_{j}-f(\underline{y},\underline{\phi}_{1},
\underline{\phi}_2){\phi}_{1j}-f(\underline{y},\underline{
\phi}_{2},\underline {\phi}_1)\phi_{2j} \bigr)
\leq\sum_{j=1}^mL(y_{j})=\sum
_{j=1}^mL\bigl(|y_{j}|\bigr)
\]
for $j=1,\ldots,m$.

We now show that
\[
\sum_{j=1}^mL\bigl(|y_{j}|\bigr)\leq L
\Biggl(\sum_{j=1}^m|y_{j}|
\Biggr).
\]
Consider $x,z\in\mathbb{R}$. Without loss of generality, we assume that
$x>z>0$. It is clear that
\[
L(x+z)-L(x)=L'(x+\lambda_1 z)z
\]
and
\[
L(z)-L(0)=L'(\lambda_2)z,
\]
where $0<\lambda_1,\lambda_2<1$. From conditions (C1) and (C2), we have
$ L(x+z)-L(x)-L(z)= [L'(x+\lambda_1z)-L'(\lambda_2z)]z \geq0$.
Thus, $\sum_{j=1}^mL(|y_{j}|)\leq L (\sum_{j=1}^m|y_{j}| )$. It
then follows that
\[
L \bigl( \bigl|y_{j}-f(\underline{y},\underline{\phi}_{1},
\underline {\phi}_2){\phi}_{1j}-f(\underline{y},\underline{
\phi}_{2},\underline {\phi}_1)\phi_{2j} \bigr|
\bigr)\leq L \Biggl(\sum_{l=1}^m|y_{l}|
\Biggr).
\]

From (C2), so we have
\[
\bigl|y_{j}-f(\underline{y},\underline{\phi}_{1},\underline{
\phi}_2){\phi}_{1j}-f(\underline{y},\underline{
\phi}_{2},\underline{\phi}_1)\phi _{2j} \bigr|\leq
\sum_{l=1}^m|y_{l}|.
\]
Furthermore, we have
\[
\bigl|f(\underline{y},\underline{\phi}_{1},\underline{\phi
}_2){\phi}_{1j}+f(\underline{y},\underline{
\phi}_{2},\underline{\phi}_1)\phi _{2j} \bigr|\leq
\sum_{l=1}^m2|y_{l}|.
\]
Also note that $\|\underline{\phi}_1\|=\|\underline{\phi}_2\|=1$ and
$\underline{\phi}_1\,\bot\,\underline{\phi}_2$, it then follows that
\begin{eqnarray*}
&&f^2(\underline{y},\underline{\phi}_{1},\underline{\phi
}_2)+f^2(\underline{y},\underline{\phi}_{2},
\underline{\phi}_1)
\\
&&\qquad =\sum_{j=1}^m \bigl[f(\underline{y},
\underline{\phi}_{1},\underline{\phi}_2){
\phi}_{1j}+f(\underline{y},\underline{\phi}_{2},\underline {
\phi}_1)\phi_{2j} \bigr]^2
\leq 4m \Biggl(\sum_{j=1}^m|y_{j}|
\Biggr)^2\leq4m^2\|\underline{y}\|^2.
\end{eqnarray*}

With the similar arguments, we have
\[
\bigl|f \bigl(\underline{y},\underline{\phi}_{1}^*,\underline{\phi
}_2^* \bigr){\phi}_{1j}^*+f \bigl(\underline{y},
\underline{ \phi}_{2}^*,\underline{\phi}_1^* \bigr)
\phi_{2j}^* \bigr|\leq\sum_{j=1}^m2|y_{j}|,
\]
where $\lambda\in(0,1)$, $\underline{\phi}_1^*=\lambda\underline
{\phi
}_1^{(1)}+(1-\lambda)\underline{\phi}_1^{(2)}$ and $\underline
{\phi
}_2^*=\lambda\underline{\phi}_2^{(1)}+(1-\lambda)\underline{\phi
}_2^{(2)}$.
Note that
\begin{eqnarray*}
&&\sum_{j=1}^m \bigl|f \bigl(\underline{y},
\underline{\phi}_{1}^*,\underline {\phi}_2^* \bigr){
\phi}_{1j}^*+f \bigl(\underline{y},\underline{\phi}_{2}^*,
\underline{\phi}_1^* \bigr)\phi_{2j}^* \bigr|^2
\\
&&\qquad = \bigl[f^2 \bigl(\underline{y},\underline{\phi}_{1}^*,
\underline{\phi}_2^* \bigr)+f^2 \bigl(\underline{y},
\underline{\phi}_{2}^*,\underline{\phi}_1^* \bigr) \bigr]
\\
&&\quad\qquad{}+2\lambda(1-\lambda)f \bigl(\underline{y},\underline{\phi}_{1}^*,
\underline {\phi}_2^* \bigr)f \bigl(\underline{y},\underline{
\phi}_{2}^*,\underline{\phi}_1^* \bigr)
\\
&&\qquad\qquad{}\times \bigl[\underline{
\phi}_1^{(1)T} \bigl(\underline{\phi}_2^{(2)}-
\underline {\phi}_2^{(1)} \bigr)+\underline{
\phi}_2^{(1)T} \bigl(\underline{\phi}_1^{(2)}-
\underline{\phi}_1^{(1)} \bigr) \bigr]
\\
&&\quad\qquad{}+2\lambda(1-\lambda)
\bigl[f^2 \bigl(\underline{y},\underline{\phi
}_{1}^*,\underline{\phi}_2^* \bigr)\underline{
\phi}_1^{(1)T} \bigl(\underline {\phi}_1^{(2)}-
\underline{\phi}_1^{(1)} \bigr)
\\
&&\hspace*{92pt}{} +f^2 \bigl(
\underline{y},\underline {\phi}_{2}^*,\underline{\phi}_1^*
\bigr)\underline{\phi}_2^{(1)T} \bigl(\underline {\phi
}_2^{(2)}-\underline{\phi}_2^{(1)}
\bigr) \bigr]
\\
&&\qquad \geq(1-3\tau/2) \bigl[f^2 \bigl(\underline{y},\underline{\phi
}_{1}^*,\underline{\phi}_2^* \bigr)+f^2 \bigl(
\underline{y},\underline{\phi}_{2}^*,\underline{\phi}_1^*
\bigr) \bigr],
\end{eqnarray*}
so it follows that
\begin{eqnarray*}
&&f^2 \bigl(\underline{y},\underline{\phi}_{1}^*,
\underline{\phi}_2^* \bigr)+f^2 \bigl(\underline{y},
\underline{\phi}_{2}^*,\underline{\phi}_1^* \bigr)
\\
&&\qquad \leq (1-3\tau/2)^{-1} \Biggl(2m\sum_{j=1}^m|y_{j}|
\Biggr)^2
\\
&&\qquad \leq 4m^3(1-3\tau/2)^{-1}\|\underline{y}\|^2.
\end{eqnarray*}\upqed
\end{pf*}

%
%
\begin{lemma}\label{Lipschitz}
If the result of Lemma~\ref{moment} and conditions \textup{(C2)--(C4)} hold, the
following inequality holds for $\underline{\phi}_1$ in a neighbor of
$\underline{\phi}_1^{(0)}$,
%
\begin{equation}
\label{Lips} \qquad\bigl|L' \bigl(y_j-f(\underline{y},
\underline{\phi}_1)\phi_{1j} \bigr)-L'
\bigl(y_j-f \bigl(\underline{y},\underline{\phi}_1^{(0)}
\bigr)\phi_{1j}^{(0)} \bigr) \bigr|\leq C\|\underline{y}\|\bigl\|
\underline{\phi}_1-\underline{\phi}_{1}^{(0)}\bigr\|
\end{equation}
for $j=1,\ldots,m$, where $y_j$ is the $j$th component of the vector
$\underline{y}$, $f$ is defined in~(\ref{implicit1})--(\ref{implicit2})
with $r=1$, and $C$ is some constant.
\end{lemma}

\begin{pf}
Without\vspace*{-2pt} loss of generality, we assume that $\phi^{(0)}_{1j}\neq0$ for
$j=1,\ldots,m_1$, and $\phi^{(0)}_{1j}=0$ for $j=m_1+1,\ldots,m$.
\begin{longlist}[(ii)]
\item[(i)]  Now we consider $j=1,\ldots,m_1$.
Consider unit vectors $\underline{\phi}$ and $\underline{\nu}$ such
that $\max\{\|\underline{\phi}-\underline{\phi}_{1}^{(0)}\|,\|
\underline
{\nu}-\underline{\phi}_1^{(0)}\|\}\leq\tau/2$,\vspace*{1pt} where $0<\tau<2/3$.

If $L'' (y_j-f(\underline{y},\underline{\phi})\phi_j )=0$, then
\[
\biggl\llvert \frac{\partial L' (y_j-f(\underline{y},\underline{\phi
})\phi
_j )}{\partial\phi_l} \biggr\rrvert =0.
\]

We\vspace*{-1pt} now consider the case where $L'' (y_j-f(\underline
{y},\underline
{\phi})\phi_j )>0$.
Let $K_1= \min\{|\phi_{1j}^{(0)}|,\break j=1,\ldots,m_1\}$.
When $\underline{\phi}$ is\vspace*{1.5pt} sufficiently close to $\underline{\phi
}{}^{(0)}_1$, we must have $|\phi_j|\geq K_1/2$, for $j=1,\ldots,m_1$. It
then follows from condition (C3) that
\[
\sum_{j=1}^mL''
\bigl(y_j-f(\underline{y},\underline{\phi})\phi _j
\bigr)\phi_j^2>0.
\]

Note that
\[
\sum
_{j=1}^mL'
\bigl(y_{j}-f(\underline{y},\underline{\phi})\phi _{j}
\bigr)\phi_{j}=0,
\]
based on the definition of $f$.
By the implicit function theorem, the partial derivatives of $f$ with
respect to $\underline{\phi}$ is
%
\begin{equation}
\label{partial} \frac{\partial f(\underline{y},\underline{\phi})}{\partial\phi
_j}=-\frac{L' (y_j-f(\underline{y},\underline{\phi})\phi_j
)-L'' (y_j-f(\underline{y},\underline{\phi})\phi_j
)f(\underline
{y},\underline{\phi})\phi_j} {
\sum_{j=1}^mL'' (y_j-f(\underline{y},\underline{\phi})\phi
_j
)\phi_j^2}
\end{equation}
for $j=1,\ldots,m$.\vadjust{\goodbreak}

Consider the partial derivative
\begin{eqnarray*}
&& \frac{\partial L' (y_j-f(\underline{y},\underline{\phi})\phi
_j)}{\partial\phi_l}
\\
&&\qquad = \cases{ \displaystyle -L''
\bigl(y_j-f(\underline{y},\underline{\phi})\phi_j \bigr)
\phi_j\frac{\partial f(\underline{y},\underline{\phi})}{\partial\phi_l},&\quad $j\neq l$,
\vspace*{3pt}\cr
\displaystyle
-L'' \bigl(y_j-f(\underline{y},\underline{
\phi})\phi_j \bigr) \biggl\{\phi_j\frac{\partial f(\underline{y},\underline{\phi})}{\partial\phi_l}+f(
\underline{y},\underline{\phi}) \biggr\},&\quad $j = l$.}
\end{eqnarray*}

Let $K_2=K_1/2$ and
\[
z_j(\underline{y},\underline{\phi})=\frac{L'' (y_j-f(\underline
{y},\underline{\phi})\phi_j )\phi_j}{\sum_{l=1}^mL''
(y_l-f(\underline{y},\underline{\phi})\phi_l )\phi_l^{2}}.
\]
Consider
%
%
\begin{eqnarray*}
\hspace*{-5pt}&&\bigl|z_j(\underline{y},\underline{\phi})\bigr|
\\
\hspace*{-5pt}&&\!\qquad =K_2^{-1}
\bigl(L'' \bigl(y_j-f(\underline{y},\underline{\phi})\phi
_j\bigr)|\phi_j|\bigr)
\\
\hspace*{-5pt}&&\!\quad\qquad{} \Bigg/\Biggl(\sum_{j=1}^{m_1}L'' \bigl(y_j-f(\underline{y},\underline
{\phi})\phi_j \bigr)\phi_j^{2}/K_2+\sum_{j=m_1+1}^{m}L''
\bigl(y_j-f(\underline{y},\underline{\phi})\phi_j \bigr)\phi_j^{2}/K_2\Biggr)
\\
\hspace*{-5pt}&&\!\qquad \leq K_2^{-1}\frac{L'' (y_j-f(\underline{y},\underline{\phi
})\phi_j )|\phi_j|}{\sum_{j=1}^{m_1}L'' (y_j-f(\underline
{y},\underline
{\phi})\phi_j )|\phi_j|}\leq K_2^{-1}.
\end{eqnarray*}

It then follows from assumption \textup{(C4)} and Lemma~\ref{moment} that
%
\begin{equation}
\label{boundderiv} \biggl\llvert \frac{\partial L' (y_j-f(\underline{y},\underline{\phi
})\phi
_j )}{\partial\phi_l} \biggr\rrvert \leq
C_1 \bigl\{\bigl|f(\underline {y},\underline {\phi})\bigr|+\|\underline{y}\|
\bigr\}\leq C\| \underline{y}\|
\end{equation}
for some constant $C$.
Hence, by (C2)--(C4), we obtain
%
\begin{equation}
\label{Lips1} \bigl|L' \bigl(y_j-f(\underline{y},
\underline{\phi})\phi_j \bigr)-L' \bigl(y_j-f(
\underline{y},\underline{\nu})\nu_j \bigr) \bigr|\leq C\| \underline {y}
\| \|\underline{\phi}-\underline{\nu}\|. 
\end{equation}

\item[(ii)] Now consider $j=m_1+1,\ldots,m$. By condition (C4), we have
\[
L''(x)\leq C_0
\]
for some constant $C_0$, and $x\in\mathbb{R}$.
It follows from condition (C3) that
\begin{eqnarray*}
&& \bigl|L' \bigl(y_j-f(\underline{y},\underline{\phi})
\phi_j \bigr)-L' \bigl(y_j-f(\underline{y},
\underline{\nu})\nu_j \bigr) \bigr|
\\
&&\qquad \leq C_0\bigl|f(\underline{y},\underline{\phi})\phi_j-f(
\underline {y},\underline{\nu})\nu_j\bigr|
\\
&&\qquad \leq C_0 \bigl\{\bigl|f(\underline{y},\underline{\phi})-f(\underline
{y},\underline {\nu})\bigr||\phi_j|+\bigl|f(\underline{y},\underline{\nu})\bigr||
\phi_j-\nu _j| \bigr\}
\\
&&\qquad \leq C_0 \bigl[ \bigl\{\bigl|f(\underline{y},\underline{\phi})\bigr|+\bigl|f(
\underline {y},\underline{\nu})\bigr| \bigr\}\bigl|\phi_j-\phi_j^{(0)}\bigr|+\bigl|f(
\underline {y},\underline {\nu})\bigr||\phi_j-\nu_j| \bigr].
\end{eqnarray*}
%
It then follows from Lemma~\ref{moment} that
%
\begin{eqnarray}
\label{Lips2}
&& \bigl|L' \bigl(y_j-f(\underline{y},
\underline{\phi})\phi_j \bigr)-L' \bigl(y_j-f(
\underline{y},\underline{\nu})\nu_j \bigr) \bigr|
\nonumber\\[-8pt]\\[-8pt]
&&\qquad\leq\| \underline {y}\|
\bigl(C_2\bigl|{\phi_j}-\phi_j^{(0)}\bigr|+C_3|
\phi_j-\nu_j| \bigr)\nonumber
\end{eqnarray}
for some constants $C_2$ and $C_3$.

Thus, by (\ref{Lips1}) and (\ref{Lips2}), we obtain (\ref{Lips}).\quad\qed
\end{longlist}\noqed
\end{pf}

\begin{pf*}{Proof of Theorem~\ref{sufficiency}}
By (\ref{Lips}) and Lemma~4.6 of \citet{HeAndShao}, we have
%
\begin{eqnarray}\label{ext3}
\quad&& \sup_{|\varphi-\varphi_0|\leq Cn^{-1/2}} \Biggl\llvert \sum
_{i=1}^{n}a_i \bigl[\gamma (
\underline{y}_i,\varphi)-\gamma(\underline{y}_i,
\varphi_0)-E \bigl\{ \gamma (\underline{y}_i,\varphi)-
\gamma(\underline{y}_i,\varphi_0) \bigr\} \bigr] \Biggr
\rrvert
\nonumber\\[-14pt]\\[-2pt]
&&\qquad  =O_p \bigl(n^{1/2} \bigr),\nonumber
\end{eqnarray}
where $\gamma(\underline{y}_i,\varphi)=\sum_{j=1}^mL'
(y_{ij}-f(\underline{y}_i,{\underline{\phi}_1}){\phi}_{1j}
){\phi}_{2j}$.

By the\vspace*{1pt} similar arguments to those used to obtain (\ref{relation1}) and
(\ref{relation2}), we obtain $f(\underline{y}_i,\underline{\phi
}_1^{(0)})=\theta_{1i}^{(0)}+f(\theta_{2i}^{(0)}\underline{\phi
}_2^{(0)}+\underline{\varepsilon}_i,\underline{\phi}_1^{(0)})$.
It\vspace*{1pt} then follows from conditions \mbox{(C1)--(C4),} (D2) and (\ref{partial}) that
%
\begin{equation}
\label{ext4} n^{-1}\sum_{i=1}^{n}a_i
\frac{\partial E\{ \gamma(\underline
{y}_i,\varphi_0)\}}{\partial\varphi_0}=n^{-1}\sum_{i=1}^{n}a_iE
\biggl\{ \frac{\partial\gamma(\underline{y}_i,\varphi_0)}{\partial\varphi
_0} \biggr\}=o(1),
\end{equation}
when $\underline{a}\,\bot\,\underline{\mu}_1$ and $\underline{\mu
}_2=\underline{0}$.
From (\ref{boundderiv}) and (C4), we know that
\[
\biggl\llvert \frac{\partial\gamma(\underline{y}_i,\varphi)}{\partial
\varphi
} \biggr\rrvert \leq C_1\|
\underline{y}\|
\]
for some constants $C_1$.
It then follows from condition (C3) and the moment condition on
$\underline{y}_i$ that
$n^{-1}\sum_{i=1}^{n}a_i\frac{\partial E\{ \gamma(\underline
{y}_i,\varphi)\}}{\partial\varphi}$ uniformly converges to a
continuous\vadjust{\goodbreak} function.
Thus, it follows from (\ref{ext3}) and (\ref{ext4}) that
%
\begin{eqnarray}
\label{approx} &&\sum_{i=1}^na_i
\Biggl\{\sum_{j=1}^mL'
\bigl(y_{ij}-f(\underline {y}_i,\tilde{\underline{
\phi}_1})\tilde{\phi}_{1j} \bigr)\tilde {\phi
}_{2j} \Biggr\}
\nonumber\\[-8pt]\\[-8pt]
&&\qquad =\sum_{i=1}^na_i
\Biggl\{\sum_{j=1}^mL'
\bigl(y_{ij}-f \bigl(\underline {y}_i,{\underline{
\phi}_1^{(0)}} \bigr){\phi}_{1j}^{(0)}
\bigr){\phi}_{2j}^{(0)} \Biggr\}+o_p
\bigl(n^{1/2} \bigr).\nonumber
\end{eqnarray}

Under condition (D2) and the null hypothesis that $\underline{\mu
}_2=\underline{0}$, we have
\begin{eqnarray*}
&&n^{-1/2}\sum_{i=1}^na_i
\gamma(\underline{y}_i,\varphi_0)
\\
&&\qquad =n^{-1/2}\sum_{i=1}^na_i
\Biggl\{\sum_{j=1}^mL' \bigl(
\theta _{2i}^{(0)}\phi_{2j}^{(0)}+
\varepsilon_{ij}-f \bigl(\theta _{2i}^{(0)}\underline{
\phi}_{2}^{(0)}+\underline{\varepsilon}_i,{
\underline{\phi}_1^{(0)}} \bigr){\phi}_{1j}^{(0)}
\bigr){\phi}_{2j}^{(0)} \Biggr\}
\\
&&\qquad \stackrel{L} {\longrightarrow} N \bigl(0,\alpha^2 \bigr)
\end{eqnarray*}
as $n\rightarrow\infty$, where
\[
\alpha^2=\operatorname{Var} \Biggl\{\sum_{j=1}^mL'
\bigl( \bigl(\theta_{2i}^{(0)}-\mu _{2i} \bigr)\phi
_{2j}^{(0)}+\varepsilon_{ij}-f \bigl( \bigl(
\theta_{2i}^{(0)}-\mu _{2i} \bigr)\underline{
\phi}_{2}^{(0)}+\underline{\varepsilon}_i,{
\underline {\phi}_1^{(0)}} \bigr){\phi}_{1j}^{(0)}
\bigr){\phi}_{2j}^{(0)} \Biggr\}.
\]

Note that
\[
\bigl|\gamma^2(\underline{y}_i,\tilde\varphi)-
\gamma^2(\underline {y}_i,\varphi_0)\bigr|\leq K
\|\underline{y}_i\|\|\tilde{\varphi}-\varphi _0\|
\]
for some constant $K$,
so $\tilde{\sigma}_n^2\stackrel{p}{\longrightarrow} \alpha^2$ as
$n\rightarrow\infty$, under the null that $\underline{\mu
}_2=\underline
{0}$. Therefore, we obtain (\ref{asym1}).
\end{pf*} 
\end{appendix}

\section*{Acknowledgements}
We thank anonymous reviewers, including an Associate Editor, who
provided valuable critiques about earlier versions of this paper. Their
comments motivated us to find better robust low-rank approximations to
data matrices with a solid theoretical underpinning. They also helped
us improve the presentation.


\begin{supplement} \label{supp}
\stitle{Additional details of case study and technical proofs\\}
\slink[doi]{10.1214/13-AOS1186SUPP} 
\sdatatype{.pdf}
\sfilename{aos1186\_supp.pdf}
\sdescription{We provide details of the case study in Section~\ref{application} and
complete the proofs of technical lemmas, as well as Theorems \ref{AsympNorm}--\ref{Consistrow}
and \ref{test2} of this paper.}
\end{supplement}

%

\printaddresses

\end{document}